\input ./style/arxiv-general.cfg
\documentclass[MSNbibl,number,citesort,dvips]{arxbj}
\makeatletter
   \@ifpackageloaded{graphicx}{}{\usepackage{graphicx}}
\makeatother
\usepackage{upgreek}

\aid{0}
\volume{21}
\issue{4}
\pubyear{2015}
\firstpage{2336}
\lastpage{2350}
\doi{10.3150/14-BEJ646} 
\docsubty{FLA}

\makeatletter

\newcommand{\nfrac}[2]{#1/#2}

\newcommand{\afrac}[2]{#1/(#2)}

\newcommand{\id}{\operatorname{Id}}
\newcommand{\eqref}[1]{(\ref{#1})}
\newcommand{\be}{\mathbf{E}}
\newcommand{\fh}{\mathfrak{H}}
\newcommand{\bp}{\mathbf{P}}
\newcommand{\var}{\operatorname{\mathbf{Var}}}
\newcommand{\cg}{{\mathcal G}}
\newcommand{\ch}{{\mathcal H}}
\newcommand{\cn}{{\mathcal N}}
\newcommand{\al}{\alpha}
\newcommand{\ep}{\varepsilon}
\newcommand{\la}{\lambda}
\newcommand{\si}{\sigma}
\newcommand{\N}{{\mathbb N}}
\newcommand{\R}{{\mathbb R}}
\newcommand{\Z}{{\mathbb Z}}
\newtheorem{theorem}{Theorem}[section]
\newtheorem{corollary}[theorem]{Corollary}
\newtheorem{hypothesis}[theorem]{Hypothesis}
\newtheorem{lemma}[theorem]{Lemma}
\newtheorem{proposition}[theorem]{Proposition}
\newremark{remark}[theorem]{Remark}
\newremark{example}[theorem]{Example}
\makeatother

\begin{document}
\begin{frontmatter}

\title{Density convergence in the Breuer--Major theorem for
Gaussian stationary sequences}
\runtitle{Limit theorems for Gaussian stationary sequences}

\begin{aug}
\author[A]{\inits{Y.}\fnms{Yaozhong} \snm{Hu}\corref{}\thanksref{A,e1}\ead[label=e1,mark]{yhu@ku.edu}},
\author[A]{\inits{D.}\fnms{David} \snm{Nualart}\thanksref{A,e2}\ead[label=e2,mark]{nualart@ku.edu}},
\author[B]{\inits{S.}\fnms{Samy} \snm{Tindel}\thanksref{B}\ead[label=e3]{samy.tindel@univ-lorraine.fr}}
\and
\author[C]{\inits{F.}\fnms{Fangjun} \snm{Xu}\thanksref{C}\ead[label=e4]{fangjunxu@gmail.com}}
\address[A]{Department of Mathematics, University of Kansas, Lawrence, KS
66045, USA.\\
\printead{e1,e2}}
\address[B]{Institut \'{E}lie Cartan, Universit\'{e} de Lorraine,
54506 Vandoeuvre-l\`{e}s-Nancy, France.\\
\printead{e3}}
\address[C]{School of Finance and Statistics, East China Normal
University, Shanghai 200241, China.\\
\printead{e4}}
\end{aug}

\received{\smonth{3} \syear{2014}}
\revised{\smonth{5} \syear{2014}}

%
\begin{abstract}
Consider a Gaussian stationary sequence with unit variance $X=\{X_k;
k\in\N\cup\{0\}\}$. Assume that the central limit theorem holds for a
weighted sum of the form $V_n=n^{-1/2}\sum^{n-1}_{k=0} f(X_k)$, where
$f$ designates a finite sum of Hermite polynomials. Then we prove
that the uniform convergence of the density of $V_{n}$ towards the
standard Gaussian density also holds true, under a mild additional
assumption involving the causal representation of $X$.
\end{abstract}

%
\end{frontmatter}

\section{Introduction}
Let $X=\{X_k; k\in\mathbb{N} \cup\{0\} \}$ be a centered Gaussian stationary
sequence with unit variance. For all $v\in\Z$, we set
$\rho(v)=\be[X_{0}X_{|v|}]$. Therefore, $\rho(0)=1$ and
$|\rho(v)|\leq1$ for all $v$. Let $\gamma$ be the standard Gaussian
probability measure and $f\in L^2(\gamma)$ be a fixed deterministic
function such that $\be[f(X_1)]=0$. We expand $f$ in the orthonormal
basis of Hermite polynomials $\{H_{k}; k\ge0\}$, which are more
specifically defined in Section~\ref{ss:mall}. In particular, if $f$
has Hermite rank $d\geq1$, it admits the following Hermite
expansion:
\[
f(x)=\sum^{\infty}_{j=d} a_j
H_j(x),
\]
with $a_d\neq0$. Define $V_n=\frac{1}{\sqrt{n}}\sum^{n-1}_{k=0}
f(X_k)$. Then the celebrated Breuer--Major theorem (see
\cite{breuermajor} or Theorem~7.2.4 in \cite{nourdinpeccatibook})
can be written as follows:

\begin{theorem}\label{thm:breuer-major}
Suppose that $\sum_{v\in\Z} |\rho(v)|^d<\infty$ and set
$\sigma^2=\sum^{\infty}_{j=d}j! a^2_j \sum_{v\in\Z} \rho(v)^j$,
which is assumed to be in $(0,\infty)$. Then the convergence
%
\begin{equation}
\label{breuermajor} V_n\stackrel{\mathit{Law}} {\longrightarrow} \cn
\bigl(0,\sigma^2\bigr)
\end{equation}
holds true as $n$ tends to infinity.
\end{theorem}

We shall be in fact interested in a particular case of Theorem~\ref{thm:breuer-major} for finite linear combinations of Hermite
polynomials, which is stated here for convenience.

\begin{corollary}\label{cor:breuer-major-Hq}
Consider $2\le d \le q<\infty$ and a family of real numbers
$\{a_{j} ; j=d,\ldots,q\}$. Let $H_{j}$ be the $j$th order
Hermite polynomial, and assume that $\si^{2}\in(0,\infty)$, where
$\si^{2}\equiv\sum^{q}_{j=d}j! a^2_j \sum_{v\in\Z} \rho(v)^j$. Set
%
\begin{equation}
\label{eq:def-Vnq} V_{n}^{d,q}=\frac{1}{\sqrt{n}}\sum
^{n-1}_{k=0}\sum_{j=d}^{q}
a_{j} H_{j}(X_k).
\end{equation}
Then $V^{d,q}_n\stackrel{\mathit{Law}}{\longrightarrow} \cn
(0,\sigma^2)$ as $n$ tends to infinity. In particular, we have
%
\begin{equation}
\label{eq:def-alpha-n} \lim_{n\to0} \be \bigl[ \bigl( V_{n}^{d,q}
\bigr)^{4} \bigr]=3\sigma^4.
\end{equation}
\end{corollary}

\begin{remark}
The relation between Gaussian convergence in law for sequences in a
fixed Wiener chaos and behavior of the fourth moment has been extensively
studied since the seminal paper~\cite{nualartpeccati}. We will need
only a small part of the information available on the topic, such as
relation~\eqref{eq:def-alpha-n}.
\end{remark}

Due to its importance, Breuer--Major theorem has been extended and
refined in several directions. Important generalizations can be found
in Arcones \cite{Ar} (multidimensional case), Chambers and Slud \cite{ChSl} and Giraitis and Surgailis \cite{GiSu}. A proof of Theorem~\ref
{thm:breuer-major} using a combination of Stein's method with Malliavin
calculus was given by Nourdin, Peccati and Podolskij in \cite{NouPePo}, where one can find explicit bounds in the total variation
and Wasserstein distances. We refer the reader to the monograph by
Nourdin and Peccati \cite{nourdinpeccati} for a more detailed account
on this topic.

We shall mainly be concerned here by
convergences of densities, and here again the relationship between
fourth moment behavior and various type of convergences of random
variables in a fixed Wiener chaos have been thoroughly studied in the
recent past. The interested reader is referred to
\cite{nourdinpeccatibook} for further details, but we will use here
the following recent criterion (see \cite{nourdinnualart},
Corollary~1.2, and \cite{hulunualart}, Corollary~4.6).

\begin{theorem}\label{thm:density-fourth-moment}
Let $\{F_{n}; n\in\N\}$ be a sequence of random variables
belonging to a fixed chaos $\ch_{q}$
with
$q\ge2$. Suppose $ \be[F_n^2]=1$ and $\lim_{n\rightarrow\infty}
\be[F_n^4] = 3 $.
Let $p_{F_n}$ be the density of
the random variable $F_n$ and let $\phi(x)=(2\uppi)^{-1/2}
\exp(-|x|^2/2)$ be the density of the standard Gaussian distribution
on $\R$.
\begin{longlist}[(ii)]
\item[(i)] Suppose that for some $\varepsilon>0$,
\[
\sup_{n } \be \bigl[\|DF_n\|^{-4-\ep}
\bigr]<\infty.
\]
Then, there exists a constant $c$ such that for all $n\ge1$,
\[
\sup_{x\in\R} \bigl|p_{F_{n}}(x) - \phi(x) \bigr| \le c \sqrt
{\be\bigl[F_n^4\bigr] - 3 } .
\]
\item[(ii)] Suppose that for all $p\ge1$,
\[
\sup_{n } \be \bigl[\|DF_n\|^{-p}
\bigr]<\infty.
\]
Then, for any $m\ge0$, there exists a constant $c_m$ such that for all
$n\ge1$,
\[
\sup_{x\in\R} \bigl|p^{(m)}_{F_{n}}(x) -
\phi^{(m)}(x) \bigr| \le c_m \sqrt {\be
\bigl[F_n^4\bigr] - 3 } .
\]
\end{longlist}
\end{theorem}

The goal of the current paper is to apply the criterion given by
Theorem~\ref{thm:density-fourth-moment} in order to get convergence
of density in the landmark of Breuer--Major theorem. In order to do
this, we need a uniform estimate on the negative moments
of the Malliavin derivative of the sequence, and this is the contents
of our main result.

\begin{theorem}\label{moment}
Let $X$ be a Gaussian stationary sequence whose spectral density
$f_\rho$ satisfies $
\log(f_{\rho})\in L^{1}([-\uppi,\uppi])$ (see Hypothesis
\ref{hyp:spectral-density} and the examples in the next section). Let
$V_{n}^{d,q}$ be the random variable defined by \eqref{eq:def-Vnq},
and assume the hypothesis of Corollary~\ref{cor:breuer-major-Hq} to
be satisfied. Then for any $p\geq1$, there exists $n_0$ such that
%
\begin{equation}
\label{eq:bnd-DV-neg-moments} \sup_{n\geq n_0} \be \bigl[\bigl\|DV_{n}^{d,q}
\bigr\|^{-p} \bigr]<\infty.
\end{equation}
\end{theorem}

In the case of a fixed Wiener chaos, we can obtain the following consequence.

\begin{corollary}
Under the conditions of Theorem~\ref{moment}, if $q=d$, and we define
$F_n= V^{d,d}_n /\sigma_n$, where $\sigma_n^2= \be[(V^{d,d}_n)^2]$,
then, for all $m\ge0$ there exists an $n_0$ (depending on $m$) such that
\[
\sup_{n\ge n_0} \sup_{x\in\R} \bigl|p^{(m)}_{F_{n}}(x)
- \phi^{(m)} (x) \bigr| \le c_m \sqrt{\be
\bigl[F_n^4\bigr] - 3 } .
\]
\end{corollary}

In the case $q\neq d$, Theorem~\ref{thm:density-fourth-moment} cannot
be applied. In the reference \cite{hulunualart} one can find results
on the uniform convergence of density for general random variables
similar to those stated in Theorem~\ref{thm:density-fourth-moment},
but they require a uniform lower bound for the negative moments of the
random variable $|\langle DF_n, -DL^{-1} F_n \rangle_\mathfrak{H} |$,
and our approach does not seem to work in this case because it is not
clear how to express $\langle DF_n, -DL^{-1} F_n \rangle_\mathfrak
{H}$ as a sum of squares. Nevertheless, condition (\ref
{eq:bnd-DV-neg-moments}) allows us to derive the uniform convergence of
the densities and their derivatives from a general result proved below
(see Proposition~\ref{prop1}) although in this case we have no
information about the rate of convergence.

\begin{corollary}
Under the conditions of Theorem~\ref{moment}, if we define
$F_n= V^{d,q}_n /\sigma_n$, where $\sigma_n^2= \be[(V^{d,q}_n)^2]$,
then, for all $m\ge0$ we have
\[
\lim_{n\rightarrow\infty} \sup_{x\in\R} \bigl|p^{(m)}_{F_{n}}(x)
- \phi^{(m)} (x) \bigr| =0.
\]
\end{corollary}

Notice that a particular case of Theorem~\ref{moment} has been
established in \cite{nourdinnualart}, for $q=2$ and
$X_{k}=B_{k+1}-B_{k}$ for a fractional Brownian motion $B$ with
Hurst parameter $H\in(0,1)$. The proof of the existence of negative
moments for $\|DF_{n}\|$ there is based on the Volterra
representation of~$B$, which leads to long computations. In
comparison our current Theorem~\ref{moment} is more general, since it
is valid for a wide class of Gaussian stationary sequences.
Its proof is also significantly simplified. These are achieved by
the introduction of two new ingredients in the proof, namely:
\begin{itemize}
\item
A general formula to compute conditional expectations for random
variables of the form $H_{q}(X_{k})$.

\item
Related to the previous item, we heavily resort to the causal
representation of $X_{k}$, which is particularly convenient in
order to compute conditional expectations.
\end{itemize}
%

Here is how our paper is structured: we give some preliminary
results concerning Gaussian stationary sequences and related
Malliavin calculus in Section~\ref{sec:preliminaries}. We then prove
our main Theorem~\ref{moment} in Section~\ref{sec:proof}.

\section{Preliminaries}\label{sec:preliminaries}

This section is devoted to some preliminaries on causal or moving
average representations for Gaussian stationary sequences, as well
as Malliavin calculus tools which will be used in the sequel.

\subsection{Moving average representation}\label{sec:moving-average-dcp}

The classical results on time series presented here are borrowed
from \cite{Be,BD}, to which we refer for further details. Start
from our Gaussian stationary sequence $\{X_{k}; k\in\N\cup\{0\}\}$ with
covariance function $\rho$. We will work under the following assumptions:

\begin{hypothesis}\label{hyp:spectral-density}
We suppose that $\rho$ admits a spectral density $f_{\rho}$
such that $\log(f_{\rho})\in  L^{1}([-\uppi,\uppi])$.
\end{hypothesis}

Condition $\log(f_{\rho})\in L^{1}([-\uppi,\uppi])$ is
referred to as \emph{purely nondeterministic} property in the
literature.
The interest of dealing with purely nondeterministic sequences is
that they admit a so-called causal representation which is
particularly convenient for conditional expectation computations.
Let us state a result in this direction, which is taken from
\cite{BD}, Chapter~5.

\begin{proposition} \label{prop1a}
Let $X$ be a Gaussian stationary sequence satisfying Hypothesis
\ref{hyp:spectral-density}. Then for each $k\in\N\cup\{0\}$ the random
variable $X_{k}$ can be decomposed as
%
\begin{equation}
\label{eq:wold-dcp} X_{k} = \sum_{j\ge0}
\psi_{j} w_{k-j},
\end{equation}
where $(w_{k})_{k\in\Z}$ is a discrete Gaussian white noise and the
coefficients $\psi_{j}$ are deterministic. We can always choose the
white noise and the coefficients in such a way that $\psi_0>0$.
\end{proposition}

\begin{pf}
The existence of the causal moving average representation is a
classical result that can be found, for instance, in
\cite{BD}, Theorem~7.5.2.
\end{pf}

Notice that from the expansion (\ref{eq:wold-dcp}) we easily obtain:
%
\begin{equation}
\label{for1} \rho(k_{1}-k_{2})=\sum
_{l=-\infty}^{k_{1}\wedge k_{2}} \psi_{k_{1}-l}
\psi_{k_{2}-l}, \quad\quad\mbox{for all } k_{1},k_{2}\in\N,
\end{equation}
and this relation will be used in the proof of our main theorem.

Let us now turn to examples for which our standing assumptions of
Hypothesis \ref{hyp:spectral-density} are
met. The following proposition
provides two typical and classical cases for which a spectral density
exists and satisfies some integrability properties.

\begin{proposition}\label{prop:relations-rho-psi}
Let $\rho$ be the covariance function of $X$. We have the following
statements.
\begin{longlist}[(ii)]
\item[(i)] If $\rho\in\ell^{1}$, then the spectral
density $f_{\rho}$ exists and is a nonnegative bounded function
defined on
$[-\uppi,\uppi]$.

 \item[(ii)]  Suppose that
$\{ \rho(k) k ^\alpha, k>0\} $ is positive and it is normalized
slowly varying at infinity for some $\alpha\in(0,1)$. That is, for
every $\delta>0$, for sufficiently large $ k $, $\rho(k) k ^{\alpha
+\delta}$ is increasing and $\rho(k) k ^{\alpha-\delta}$ is
decreasing. Then the spectral density exists and
satisfies
$\lim_{\la\to0}|\la|^{1-\al}f_{\rho}(\la)=c_{f}$ for some
constant $c_f>0$ (see \cite{gu}).
\end{longlist}
\end{proposition}

We now give two specific and important examples which satisfy
Hypothesis~\ref{hyp:spectral-density}.

\begin{example}
The so-called autoregressive
fractionally integrated moving-average (ARFIMA) processes are introduced
in \cite{granger} and \cite{hosking}. Denote by $B$ the one lag
backward operator ($BX_k=X_{k-1}$). Let $\phi(z)$ and $\theta(z)$ be
two polynomials which have no common zeros and such that the zeros of
$\phi$
lie outside the closed unit disk $\{z , |z|\le1\}$. Suppose that $X_k$ is
given by
%
\begin{equation}
\label{eq:def-arfima} \phi(B) X_k=(\id-B)^{-d} \theta(B)
w_k ,
\end{equation}
where $-1<d<1/2$, and where the operator $(\id-B)^{-d}$ is defined by:
\[
(\id-B)^{-d}=\sum_{j=1}^\infty
\eta_j B^j\quad\quad\mbox{with } \eta_j=\frac{\Gamma(d+j)}{\Gamma(j+1)\Gamma(d)}
.
\]
Also notice that in \eqref{eq:def-arfima}
the sequence $(w_k) _{k\in\mathbb{Z}}$ is a discrete Gaussian white noise.
It is well-known (see \cite{palma}, Theorem~3.4 and equation (3.19))
that under the above conditions, $\{X_k, k\in\N\}$ admits a spectral
density whose exact expression is:
\[
f(\la)=\frac{1}{2\uppi} \biggl[ 2\sin\frac{\lambda}{2} \biggr]^{-2d}
\frac{|\theta(\mathrm{e}^{-\mathrm{i}\lambda})|^2}{
|\phi(\mathrm{e}^{-\mathrm{i}\lambda})|^2}.
\]
It is thus readily checked that
Hypothesis \ref{hyp:spectral-density} is satisfied, and hence $X_{k}$
has a causal representation.
\end{example}

\begin{example}
Our second example is the fractional Gaussian noise. Let $\{B_t,
t\ge0\}$ be a fractional Brownian motion of Hurst parameter $H\in
(0, 1)$. Then $\{X_k=B_{k+1}-B_k , k\in\N\cup\{0\}\}$ is a stationary
Gaussian process with correlation
\[
\rho(k)=\frac{1}2 \bigl[ (k+1)^{2H} -2k^{2H}+(k-1)^{2H}
\bigr] .
\]
Its spectral density (see, e.g., \cite{Be}, equation (2.17)) is:
\[
f(\la)=\frac{1}{2\uppi} \sum_{k=-\infty}^\infty
\rho\bigl(|k|\bigr) \mathrm{e}^{\mathrm{i}k\lambda}=2c_f \bigl(1-\cos(\lambda)\bigr)
\sum_{j=-\infty}^\infty|2\uppi j+
\lambda|^{-2H-1} , \quad\quad\lambda\in[-\uppi, \uppi] ,
\]
where $c_f=(2\uppi)^{-1} \sin(\uppi H)\Gamma(2H+1)$. If $H\le1/2$,
it is clear that $\sum_{k=-\infty}^\infty|\rho(|k|)|<\infty$. This
implies
\[
\sup_{\lambda\in[-\uppi, \uppi]} \bigl|f(\la)\bigr|<\infty.
\]
If $1/2<H<1$, then
\[
0\le f(\la)\le 2c_f \bigl(1-\cos(\lambda)\bigr)|
\la|^{-2H-1} +2c_f \sum_{j\neq0} |2
\uppi j+\lambda|^{-2H-1} , \quad\quad\lambda\in[-\uppi, \uppi] .
\]
The first term is in $L^1 $ since $H<1$. When $j\neq0$,
$\int_{-\uppi}^\uppi|2\uppi j+\lambda|^{-2H-1} \,\mathrm{d}\lambda
\le C j^{-2H}$
for some positive constant $C$. Thus $\int_{-\uppi}^\uppi\sum_{j\neq0}
|2\uppi j+\lambda|^{-2H-1} \,\mathrm{d}\lambda<\infty$, owing to
the fact that $H>1/2$. Therefore, we have
$f\in L^1$. Summarizing we have
$f\in L^1$ for all $H\in(0, 1)$. This also implies
$\log^+ f(\lambda) \in L^1$. To see $\log^- f(\lambda)\in
L^1$, we notice that
\[
f(\la) \ge2c_f \bigl(1-\cos(\lambda)\bigr)|\la|^{-2H-1} .
\]
So $\log^- f(\lambda) \le C+ |\log [(1-\cos
(\lambda))|\la|^{-2H-1} ] |$ which is in $L^1$. In
conclusion, the sequence $X$ satisfies Hypothesis \ref{hyp:spectral-density}.
\end{example}

\subsection{Malliavin calculus}\label{ss:mall}

We start by briefly recalling some basic notation and results
connected to Gaussian analysis and Malliavin calculus. The reader
is referred to \cite{nourdinpeccatibook,nualartbook} for details or
missing proofs.

\subsubsection{Wiener space and generalities}\label{sec:wiener space}

Let $\mathfrak{H}$ be a real separable Hilbert space with inner
product $\langle\cdot,\cdot\rangle_\mathfrak{H}$. The norm of
$\mathfrak{H}$ will be denoted by $\|\cdot\|=\| \cdot
\|_{\mathfrak{H}}$. Recall that we call {\it isonormal Gaussian
process} over $\mathfrak{H}$ any centered Gaussian family $W =
\{W(h) : h\in\mathfrak{H}\}$, defined on a probability space
$(\Omega, \mathcal{F}, \mathbf{P})$ and such that $\be[W(h) W(g)] =
\langle h, g\rangle_\mathfrak{H}$ for every $h,g\in\mathfrak{H}$.

In our application the underlying Gaussian family will be a discrete
Gaussian white noise $(w_k)_{k\in\mathbb{Z}}$. The space $\fh$ is given
here by $\fh=\ell^{2}(\Z)$ (the space of square integrable sequences
indexed by $\Z$) equipped with its natural inner product. Set
$\{\ep^{j}; j\in\Z\}$ for the canonical basis of $\ell^{2}(\Z)$,
that is $\ep^{j}_{k}=\delta_{j}(k)$.
We thus identify $w_{j}$ with
$W(\ep^{j})$. Assume from now on that our underlying $\si$-algebra
$\mathcal{F}$ is generated by $W$.

For any integer $q\in\N\cup\{0\}$, we denote by $\mathcal{H}_q$ the
$q$th {\it Wiener chaos} of $W$. We recall that $\mathcal{H}_0$ is
simply $\R$ whereas, for any $q\geq1$, $\mathcal{H}_q$ is the
closed linear subspace of $L^2(\Omega)$ generated by the family of
random variables $\{H_q(W(h)), h\in\mathfrak{H},
\|h\|_{\mathfrak{H}}=1\}$, with $H_q$ the $q$th Hermite polynomial
given by
%
\begin{equation}
\label{eq:def-hermite-polynomial} H_q(x)= (-1)^q \mathrm{e}^{\nfrac{x^2}2}
\frac{\mathrm{d}^q}{
\mathrm{d}x^q} \bigl( \mathrm{e}^{-\nfrac{x^2}2} \bigr).
\end{equation}

Let $\mathcal{S}$ be the set of all cylindrical random variables of
the form
\[
F=g\bigl(W(h_1), \dots, W(h_n)\bigr),
\]
where $n\geq1$, $h_i \in\mathfrak{H}$, and $g$ is infinitely
differentiable such that all its partial derivatives have polynomial
growth. The Malliavin derivative of $F$ is the element of
$L^2(\Omega;\mathfrak{H})$ defined by
\[
DF= \sum_{i=1}^n \frac{\partial g}{\partial x_i}
\bigl(W(h_1), \dots, W(h_n)\bigr) h_i.
\]
By iteration, for every $m\geq2$, we define the $m$th derivative
$D^mF$. This is an element of $L^2(\Omega; \mathfrak{H}^{\odot
m})$, where $\mathfrak{H}^{\odot m}$ designates the symmetric $m$th
tensor product of $\fh$. For $m\geq1$ and $p\geq1$,
$\mathbb{D}^{m,p}$ denote the closure of $\mathcal{S}$ with respect
to the norm $\| \cdot\|_{m,p}$ defined by
\[
\|F\|^p_{m,p} = \be\bigl[|F|^p\bigr] + \sum
_{j=1}^m \be \bigl[ \bigl\|D^jF\bigr\|
^p_{\mathfrak{H}^{\otimes j}} \bigr].
\]
Set $\mathbb{D}^\infty= \bigcap_{m,p} \mathbb{D}^{m,p}$.
One can then extend the definition of $D^m$ to $\mathbb{D}^{m,p}$.
When $m=1$, one simply write $D$ instead of $D^1$. As a consequence
of the hypercontractivity property of the Ornstein--Uhlenbeck
semigroup (see, e.g., \cite{nourdinpeccatibook}, Theorem~2.7.2), all
the $\| \cdot\|_{m,p}$-norms are equivalent in any {\it finite} sum
of Wiener chaoses.

Finally, let us recall that the Malliavin derivative $D$ satisfies the
following \textsl{chain
rule}: if $\varphi\dvtx \R^n\rightarrow\R$ is in $\mathcal{C}^1_b$
(i.e., belongs to the set of continuously differentiable
functions with a bounded derivative) and if
$\{F_i\}_{i=1,\ldots,n}$ is a vector of elements of $\mathbb{D}^{1,2}$,
then $\varphi(F_1,\ldots,F_n)\in\mathbb{D}^{1,2}$ and
%
\begin{equation}
\label{eq:chain-rule} D\varphi(F_1,\ldots,F_n)=\sum
_{i=1}^n \frac{\partial\varphi}{\partial x_i} (F_1,\ldots,
F_n) DF_i.
\end{equation}

\subsection{Convergence of densities}
Suppose that $F$ is a random variable in $\mathbb{D}^{\infty}$ such
that $\be[ \|DF \|^{-p} ]<\infty$ for all $p\ge1$. Then, we know that
$F$ has an infinitely differentiable density and there are explicit
formulas for the density and its derivatives (see \cite{nualartbook}, Proposition~2.1.5).
Using this result, we can establish the following criteria for
convergence of densities for random variables in a finite sum of Wiener chaoses.

\begin{proposition} \label{prop1}
Let $\{F_n; n\in\mathbb{N}\}$ be a sequence of random variables
belonging to a finite sum of Wiener chaoses $\bigoplus_{k=1}^M \mathcal
{H}_k$, which converges in law to a nonzero random variable $F_\infty
$. Suppose that, for all $p\ge1$
%
\begin{equation}
\label{c1} \sup_n \be\bigl[ \|DF_n
\|^{-p} \bigr]<\infty.
\end{equation}
Then, for all $m\ge0$ the derivative $p^{(m)}_n $ of the density of
$F_n$, converges uniformly and in $L^p(\mathbb{R})$ for all $p\ge1$
to the corresponding derivative of the density of $F_\infty$.
\end{proposition}

\begin{pf}
First, notice that by condition (\ref{c1}), the random variable $F_n$
has an infinitely differentiable density $p_n$, whose derivatives can
be expressed by
%
\begin{equation}
\label{c2} p _n ^{(m)}(x)=\be \bigl[
\mathbf{1}_{\{F>x\}} G_n^{(m)} \bigr],
\end{equation}
where the random variables $G_n^{(m)} $ are defined recursively by
$G_n^{(0)}=\delta ( \frac{DF_n}{\|DF_n\|^2_H}  )$ and
\[
G_n^{(m)} = -\delta \biggl( \frac{ G_n^{(m-1)}DF_n}{\|DF_n\|^2_H} \biggr),
\]
for any $m\ge1$.
From this formula, it follows that the derivatives $p_n^{(m)}$ are
uniformly bounded and also uniformly bounded in $L^p(\mathbb{R})$ for
all $p\ge1$. Indeed, by \cite{NouPo}, Lemma~2.4, we have $\sup_n \be
[|F_n|^r] <\infty$ for all $r\ge1$. This uniform bound on the
moments, together with the equivalence of the $\| \cdot\|_{m,p}$ norms
in any finite sum of Wiener chaoses and condition (\ref{c1}) imply
that $\sup_n \| G_n^{(m)} \|_{L^p(\Omega)}=c_{m,p} <\infty$ for all
$m\ge0$. Then, we can write from (\ref{c2})
%
\begin{equation}
\label{w1} \sup_n \sup_x
\bigl|p_n ^{(m)}(x)\bigr|\le\sup_n \be \bigl[ \bigl|
G_n^{(m)}\bigr| \bigr] =c_{m,1} <\infty,
\end{equation}
and, using the fact that $\be[G_n^{(m)}] =0$, we get:
%
\begin{eqnarray}\label{w2}
\nonumber
\sup_n \bigl|p_n ^{(m)}(x)\bigr| &\le&
\sup_n \bigl( \bp\bigl(|F_n|>|x|\bigr) \be \bigl[ \bigl|
G_n^{(m)}\bigr| ^2 \bigr] \bigr)^{\nfrac{1}2}
\\
 &\le& c_{m,2} \sup_n \sqrt{\bp\bigl(
|F_n| >|x|\bigr)} \\&\le& c_{m,2} \sup_n \be
\bigl[|F_n| ^q\bigr] ^{\nfrac{1}{2}} |x| ^{-\nfrac{q} 2}
\le c |x| ^{-\nfrac{q} 2}\nonumber,
\end{eqnarray}
for all $q\ge1$ and for some constant $c$ depending on $q$ and $m$.

By \cite{NouPo}, Theorem~3.1, the sequence $F_n$ converges in total
variation to $F_\infty$, that is, the densities $p_n$ converge in
$L^1(\mathbb{R})$ to the density $p_\infty$ of $F_\infty$. The
boundedness in $L^p(\mathbb{R})$ and the uniform bound of $p_n$ imply
that this convergence holds in $L^p(\mathbb{R})$ for any $p\ge1$.

On the other hand, the estimates (\ref{w1}) and (\ref{w2}) imply that
for any $m\ge1$ and any $p\ge1$, the sequence $p^{(m)}_n$ is
uniformly bounded in $L^p(\mathbb{R})$. Therefore, for any $m\ge1$
and any $p\ge1$ the sequence $p^{(m)}_n$ is relatively compact in
$L^p(\mathbb{R})$. Suppose that a subsequence $\{p^{(m)}_{n_k}, k\ge
1\}$ converges in $L^p(\mathbb{R})$ to some limit $\tilde
p^{(m)}_\infty$. This limit must coincide with the $m$th derivative
(in the distribution sense) of $p_\infty$, and, therefore, it is
unique. This implies that for any $m\ge1$ and any $p\ge1$ the
sequence $p^{(m)}_n$ converges in $L^p(\mathbb{R})$ to the $m$th
derivative of $p_\infty$.

Finally, the uniform convergence is also easy to establish from the
convergence of the densities in $L^p(\mathbb{R})$ for all $p\ge1$.
\end{pf}

\subsection{A key lemma}

Our future computations will heavily rely on an efficient way to
compute conditional expectations. Towards this aim, we state here
some general results. Let us start with a decomposition for Hermite
polynomials.

\begin{lemma} \label{lem:hu-polynomial} For any $q\geq1$, let $H_{q}$
be the polynomial
defined by relation \eqref{eq:def-hermite-polynomial}. Consider
$y,z\in\R$ and two real parameters $a,b$ with $a^2+b^2=1$. Then the
following relation holds true:
%
\begin{equation}
\label{hu1} H_q(ay+bz)=\sum^q_{\ell=0}
\pmatrix{q\cr\ell}  a^{q-\ell} b^{\ell} H_{q-\ell}(y)
H_{\ell}(z).
\end{equation}
\end{lemma}

\begin{pf} By the definition of the Hermite polynomials, we have
%
\begin{equation}
\label{eq:identity-exp-hermite1} \mathrm{e}^{aty-\nfrac{(at)^2}{2}}
=\sum^{\infty}_{i=0}
\frac{(at)^i }{i!} H_i(y) \quad \mbox{and}\quad \mathrm{e}^{tbz-\nfrac{(bt)^2}{2}}=\sum
^{\infty}_{j=0}\frac
{(bt)^j}{j!}
H_j(z).
\end{equation}
In the same way, we also obtain
%
\begin{equation}
\label{eq:identity-exp-hermite2} \mathrm{e}^{t(ay+bz)-t^2/2}=\sum^{\infty}_{q=0}
\frac{t^q}{q!} H_q(ay+bz).
\end{equation}

Since $a^2+b^2=1$, we obviously have
$\mathrm{e}^{aty-\nfrac{(at)^2}{2}}\mathrm{e}^{tbz-\nfrac
{(bt)^2}{2}}=\mathrm{e}^{t(ay+bz)-t^2/2}$.
Thus, multiplying the right-hand sides of both identities in
\eqref{eq:identity-exp-hermite1} we recover the right-hand side of
\eqref{eq:identity-exp-hermite2}, namely:
\[
\sum^{\infty}_{q=0}\frac{ t^q}{q!}
H_q(ay+bz)=\sum^{\infty
}_{i=0}
\frac{(at)^i}{i!}H_i(y)\sum^{\infty}_{j=0}
\frac{(bt)^j}{j!} H_j(z),
\]
which easily yields the desired identity (\ref{hu1}).
\end{pf}

With this preliminary result in hand, we are ready to state our result
on conditional expectations.

\begin{proposition} \label{hu}
Let $Y$ and $Z$ be two centered Gaussian random variables such that
$Y$ is measurable with respect to a $\si$-algebra
$\cg\subset\mathcal{F}$ and $Z$ is independent of $\mathcal{G}$.
Assume that $\be[Y^2]=\be[Z^2]=1$. Then for any $q\geq1$, and real
parameters $a,b$ such that $a^2+b^2=1$, we have:
%
\begin{equation}
\label{hu2} \be\bigl[H_q(aY+bZ)|\mathcal{G}\bigr]=a^q
H_q(Y).
\end{equation}
\end{proposition}

\begin{pf}
Apply identity \eqref{hu1} in order to decompose $H_q(aY+bZ)$. Then
identity (\ref{hu2}) follows easily from the fact that $Y$ is
$\cg$-measurable, $Z$ is independent from $\cg$ and Hermite
polynomials have 0 mean under a centered Gaussian measure except for
$H_{0}\equiv1$.
\end{pf}

\subsection{Carbery--Wright inequality}

In the proof of Theorem~\ref{moment}, we will make use of the
following inequality due to Carbery and Wright \cite{carberywright}, Theorem~8, which is recalled here for convenience.

\begin{proposition}
Let $X=(X_1,\ldots,X_n)$ be a Gaussian random vector in $\R^{n}$
and $Q\dvtx \mathbb{R}^n \rightarrow\mathbb{R}$ a polynomial of degree
at most $m$. Then there is a universal constant $c>0$ such that:
%
\begin{equation}
\label{ff1} \be\bigl[\bigl|Q(X_1, \dots, X_n)\bigr|\bigr]
^{\nfrac{1}{m}} \bp\bigl(\bigl|Q(X_1, \dots, X_n)\bigr| \leq x
\bigr) \leq c m x^{\nfrac{1}m}, \qquad\mbox{for all } x>0.
\end{equation}
\end{proposition}

\section{Proof of Theorem \texorpdfstring{\protect\ref{moment}}{1.5}}\label{sec:proof}
In this section, we will prove our main result, which amounts to
show the inequality~\eqref{eq:bnd-DV-neg-moments}. This will be
done into several steps.

\textit{Step} 1: \textit{Computation of the Malliavin norm.} Recall
that $V_{n}^{d,q}$ is defined by relation \eqref{eq:def-Vnq}, and
set for the moment $f=\sum_{j=d}^{q} a_{j} H_{j}$. Invoking
relation \eqref{eq:chain-rule}, plus the fact that $Dw_{j}=\ep^{j}$
with the notation of Section~\ref{sec:wiener space}, we get:
%
\begin{equation}
\label{eq:express-DV-n-q} DV_{n}^{d,q}=\frac{1}{\sqrt{n}}\sum
_{k=0}^{n-1} f'(X_{k}) \biggl(
\sum_{j\ge0} \psi_{j} \ep^{k-j}
\biggr) =\frac{1}{\sqrt{n}} \sum_{l=-\infty}^{n-1}
\Biggl( \sum_{k=l^{+}}^{n-1} \psi_{k-l}
f'(X_{k}) \Biggr) \ep^{l} ,
\end{equation}
where $l^{+}=\max\{{l, 0\}}$. Invoking relation (\ref{for1}), it is
thus readily
checked that:
\[
\bigl\langle DV_{n}^{d,q}, DV_{n}^{d,q}
\bigr\rangle_{\fh} =\frac{1}{n}\sum^{n-1}_{k_{1},k_{2}=0}
f'(X_{k_{1}}) \rho (k_{1}-k_{2})
f'(X_{k_{2}}) ,
\]
where we recall that $\rho$ is the covariance function of the Gaussian
stationary sequence $\{X_k; k\geq0\}$. This is consistent with
the expression found in \cite{nourdinpeccatibook}, Chapter~5.
However, in order to write the above expression as sum of some squares,
we will start directly from expression \eqref{eq:express-DV-n-q}.
Since $\{\ep^{l} ; l\in\Z\}$ is an orthonormal basis of
$\ell^{2}(\Z)$ we obtain:
\[
\bigl\langle DV_{n}^{d,q}, DV_{n}^{d,q}
\bigr\rangle_{\fh} = \frac{1}{n} \sum
^{n-1}_{\ell=-\infty} \Biggl( \sum
^{n-1}_{k=\ell^{+}} \psi _{k-\ell}
f'(X_k) \Biggr)^{2} .
\]
Rearranging terms (namely, change $k-\ell$ to $ k$ and then $n-\ell
-1$ to
$m$), we end up with:
\begin{eqnarray*}
\bigl\langle DV_{n}^{d,q}, DV_{n}^{d,q}
\bigr\rangle_{\fh} &\ge& \frac{1}{n}\sum
^{n-1}_{\ell=0} \Biggl(\sum^{n-\ell-1}_{k=0}
f'(X_{\ell+k}) \psi_{k} \Biggr)^2
\\
&=& \frac{1}{n}\sum^{n-1}_{m=0}
\Biggl(\sum^{m}_{k=0} f'(X_{n-1-(m-k)})
\psi_{k} \Biggr)^2 \equiv A_{n} .
\end{eqnarray*}
As a last preliminary step we resort to the fact that $X=\{X_k;
k\in\N\cup\{0\}\}$ is a Gaussian stationary sequence, which allows
to assert
that $A_{n}$ is identical in law to $B_{n}$ with
\[
B_{n} := \frac{1}{n}\sum^{n-1}_{m=0}
\Biggl(\sum^{m}_{k=0} f'(X_{m-k})
\psi_{k} \Biggr)^2=\frac{1}{n}\sum
^{n-1}_{m=0} \Biggl(\sum^{m}_{k=0}
f'(X_{k}) \psi_{m-k} \Biggr)^2.
\]
We will now bound the negative moments of $B_{n}$.

\textit{Step} 2: \textit{Block decomposition}. We now wish to apply
the Carbery--Wright inequality \eqref{ff1} in order to get bounds
for negative moments of $B_{n}$. However, relation \eqref{ff1} only
applies to moments of small order, and this is why we proceed to a
decomposition of $B_{n}$ into smaller blocks.

Fix thus an integer $N\geq1$ and let $M=[n/N]$ be the integer part
of $n/N$. Then $n\geq NM$ and as a consequence,
\[
B_n=\frac{1}{n}\sum^{n-1}_{m=0}
\Biggl(\sum^{m}_{k=0} f'(X_{k})
\psi_{m-k} \Biggr)^2\geq \frac{1}{n}\sum
^{N-1}_{i=0}\sum^{(i+1)M-1}_{m=iM}
\Biggl(\sum^{m}_{k=0} f'(X_{k})
\psi_{m-k} \Biggr)^2.
\]
For $i=0,\dots, N-1$, define
\[
B^i_n=\frac{1}{n}\sum
^{(i+1)M-1}_{m=iM} \Biggl(\sum^{m}_{k=0}
f'(X_{k}) \psi_{m-k} \Biggr)^2
\]
so that $ B_n\ge\sum_{i=0}^{N-1} B_{n}^i$. Then it
is readily checked that:
%
\begin{equation}
\label{ineq} (B_{n})^{-\nfrac{p}{2}} \le\prod
_{i=0}^{N-1} \bigl(B_{n}^i
\bigr)^{-\afrac{p}{2N}}.
\end{equation}
Recall once again the representation of the sequence $X$ in
\eqref{eq:wold-dcp}, and denote by $\mathcal{F}_k$ the filtration
generated by $\{w_{\ell}: \ell<k\}$. Then starting from \eqref{ineq}
we obtain:
%
\begin{eqnarray}
\label{ineq1} \be \bigl[(B_{n})^{-\nfrac{p}{2}} \bigr] &\le&\be \Biggl[
\prod_{i=0}^{N-1} \bigl(B_{n}^i
\bigr)^{-\afrac{p}{2N}} \Biggr]\nonumber
\\[-8pt]\\[-8pt]
& =& \be \Biggl[ \be \bigl[ \bigl(B_{n}^{N-1}
\bigr)^{-\afrac{p}{2N}}| \mathcal {F}_{(N-1)M} \bigr] \prod
_{i=0}^{N-2} \bigl(B_{n}^i
\bigr)^{-\afrac{p}{2N}} \Biggr].\nonumber
\end{eqnarray}

\textit{Step} 3: \textit{Application of Carbery--Wright.} Let us go
back to the particular situation of $f=\sum_{j=d}^{q} a_{j}
H_{j}$, which means in particular that $f'= \sum_{j=d}^{q} j a_{j}
H_{j-1}$. We are now in a position to apply a conditional version
of inequality \eqref{ff1} to the block
$(B_{n}^{N-1})^{-\afrac{p}{2N}}$ in \eqref{ineq1}. First, we
notice
%
\begin{equation}
\label{ineq2} \be \bigl[ \bigl(B_{n}^{N-1}
\bigr)^{-\afrac{p}{2N}}| \mathcal{F}_{(N-1)M} \bigr] \le1 +\frac{p}{2N}
\int_0^1 \bp \bigl( B_{n}^{N-1}
\leq x | \mathcal{F}_{(N-1)M} \bigr) x^{-\afrac{p}{2N} -1} \,\mathrm{d}x .
\end{equation}
Since $B_{n}^{N-1}$ is a polynomial of order $m=2(q-1)$,
Carbery--Wright's inequality (\ref{ff1}) yields:
%
\begin{equation}
\label{a} \bp \bigl( B_{n}^{N-1} \leq x |
\mathcal{F}_{(N-1)M} \bigr) \le \frac{c {
x^{\afrac{1}{2(q-1)}}}}{
  [ \be ( B_{n}^{N-1} | \mathcal{F}_{(N-1)M}  )
 ]^{\afrac{1}{2(q-1)}} }.
\end{equation}

\textit{Step} 4: \textit{Estimates for the conditional expectation.}
We now estimate the conditional expectation $\be[ B_{n}^{N-1} |
\mathcal{F}_{(N-1)M} ]$. We have:
%
\begin{eqnarray}
\be \bigl[ B_{n}^{N-1} | \mathcal{F}_{(N-1)M} \bigr]
&=&\frac{1}{n}\sum_{m=(N-1)M}^{NM-1} \be
\Biggl[ \Biggl(\sum^{m}_{k=0}
f'(X_{k}) \psi_{m-k} \Biggr)^2 \Big|
\mathcal{F}_{(N-1)M} \Biggr]
\nonumber
\\[-8pt]\\[-8pt]
&\geq&\frac{1}{n}\sum_{m=(N-1)M}^{NM-1}
A_{m},\nonumber \label{eq:low-bnd-Bn-Am}
\end{eqnarray}
where we have set
\[
A_{m}=\var \Biggl( \sum^{m}_{k=(N-1)M}
f'(X_{k}) \psi_{m-k} \Big| \mathcal{F}_{(N-1)M}
\Biggr).
\]
Furthermore, notice that
\[
f'(X_k) =f' \Biggl(\sum
^k_{\ell=-\infty} \psi_{k-i} w_i
\Biggr)=f'(Y_k+Z_k),
\]
where $Y_k=\sum^{(N-1)M-1}_{i=-\infty} \psi_{k-i} w_i$ is
$\mathcal{F}_{(N-1)M}$-measurable and $Z_k=\sum^k_{i=(N-1)M}
\psi_{k-i} w_i$ is independent of $\mathcal{F}_{(N-1)M}$.
Recalling that $f'= \sum_{j=d}^{q} j a_{j} H_{j-1}$, we can
thus resort to Lemmas \ref{lem:hu-polynomial} and~\ref{hu}. This
gives:
\[
f'(X_k)-\be\bigl[f'(X_k)|
\mathcal{F}_{(N-1)M}\bigr] =\sum^{q}_{j=d}
\sum^{j-1}_{\ell=1} j a_j \pmatrix{j-1\cr
\ell }   \sigma^{j-1-\ell}_{Y_k}H_{j-1-\ell}(
\widetilde{Y}_k)\sigma^{\ell}_{Z_k}
H_{\ell}(\widetilde{Z}_k),
\]
where $\sigma_{Y_k}=[\var(Y_k)]^{1/2}$,
$\sigma_{Z_k}=[\var(Z_k)]^{1/2}$, $\widetilde{Y}_k=Y_k/\sigma_{Y_k}$
and $\widetilde{Z}_k=Z_k/\sigma_{Z_k}$. Therefore,
\begin{eqnarray*}
A_{m} &=&\be \Biggl[ \Biggl(\sum^{m}_{k=(N-1)M}
\sum^{q}_{j=d} \sum
^{j-1}_{\ell=1} a_{j,\ell,k} H_{j-1-\ell}(
\widetilde{Y}_k) H_{\ell}(\widetilde{Z}_k)
\psi_{m-k} \Biggr)^2 \Big| \mathcal{F}_{(N-1)M} \Biggr]
\\
&=& \be \Biggl[ \Biggl(\sum^{q-1}_{\ell=1} \sum
^{m}_{k=(N-1)M} \sum
^{q}_{j=(\ell+1)\vee d} a_{j,\ell,k} H_{j-1-\ell}(
\widetilde{Y}_k) H_{\ell}(\widetilde{Z}_k)
\psi_{m-k} \Biggr)^2 \Big| \mathcal{F}_{(N-1)M} \Biggr],
\end{eqnarray*}
where we have set $a_{j,\ell,k}=j a_j \bigl({j-1\atop\ell} \bigr)
\sigma^{j-1-\ell}_{Y_k}\sigma^{\ell}_{Z_k}$.

Recall that the random variables $\widetilde{Y}_k$ are
$\mathcal{F}_{(N-1)M}$-measurable while the random variables
$\widetilde{Z}_k$ are independent of $\mathcal{F}_{(N-1)M}$. By
decorrelation properties of Hermite polynomials we thus get:
\[
A_{m}= \sum^{q-1}_{\ell=1} \be
\Biggl[ \Biggl( \sum^{m}_{k=(N-1)M} \sum
^{q}_{j=(\ell+1)\vee d} a_{j,\ell,k}
H_{j-1-\ell}(\widetilde{Y}_k) H_{\ell}(
\widetilde{Z}_k) \psi_{m-k} \Biggr)^2 \Big|
\mathcal{F}_{(N-1)M} \Biggr]
\]
and we trivially lower bound this quantity by taking the term
corresponding to $\ell=q-1$. In this situation, the sum
$\sum^{q}_{j=(\ell+1)\vee d}$ is reduced to the term corresponding
to $j=q$, and since $a_{q,q-1,k}=q a_{q} \si_{Z_{k}}^{q-1}$ we
obtain:
\begin{eqnarray*}
A_{m} &\ge&\be \Biggl[ \Biggl(\sum^{m}_{k=(N-1)M}
q a_q \sigma^{q-1}_{Z_k} H_{q-1}(
\widetilde{Z}_k) \psi_{m-k} \Biggr)^2 \Big|
\mathcal{F}_{(N-1)M} \Biggr]
\\
&=&q^2 a_q^{2} \be \Biggl[ \Biggl(\sum
^{m}_{k=(N-1)M} \sigma^{q-1}_{Z_k}
H_{q-1}(\widetilde{Z}_k) \psi_{m-k}
\Biggr)^2 \Biggr].
\end{eqnarray*}
We now invoke the identity $\be[H_{p}(\tilde{Z}_{k_{1}})
H_{p}(\tilde{Z}_{k_{2}})]= p!
(\be[\tilde{Z}_{k_{1}}\tilde{Z}_{k_{2}}])^{p}$ in order to obtain
\[
A_{m}\ge q q! a_q^{2} \sum
^{m}_{k_{1},k_{2}=(N-1)M} \sigma^{q-1}_{Z_{k_{1}}}
\sigma^{q-1}_{Z_{k_{2}}} \be [\widetilde{Z}_{k_{1}}
\widetilde{Z}_{k_{2}} ]^{q-1} \psi _{m-k_1}
\psi_{m-k_2}.
\]
Furthermore, similarly to (\ref{for1}), it is readily checked that:
\[
\be [\widetilde{Z}_{k_{1}} \widetilde{Z}_{k_{2}} ] =
\frac{1}{\sigma_{Z_{k_{1}}} \sigma_{Z_{k_{2}}}} \sum_{i=(N-1)M}^{k_{1}\wedge k_{2}}
\psi_{k_{1}-i} \psi_{k_{2}-i} ,
\]
and thus
\begin{eqnarray*}
A_{m} &\ge& q q! a_q^{2} \sum
^{m}_{k_{1},k_{2}=(N-1)M} \Biggl(\sum_{i=(N-1)M}^{k_{1}\wedge k_{2}}
\psi_{k_{1}-i} \psi _{k_{2}-i} \Biggr)^{q-1}
\psi_{m-k_1}\psi_{m-k_2}
\\
&=& q q! a_q^{2} \sum_{i_{1},\dots,i_{q-1}=(N-1)M}^{m}
\sum_{k_{1},k_{2}=\max(i_{1},\dots,i_{q-1})}^{m} \psi_{m-k_1}
\psi_{m-k_2}\prod_{j=1}^{q-1}
\psi_{k_{1}-i_{j}} \psi _{k_{2}-i_{j}}
\\
&=& q q! a_q^{2} \sum_{i_{1},\dots,i_{q-1}=(N-1)M}^{m}
\Biggl( \sum^{m}_{k=\max(i_{1},\dots,i_{q-1})}
\psi_{m-k}\prod_{j=1}^{q-1}
\psi_{k-i_{j}} \Biggr)^{2}.
\end{eqnarray*}
Here again, this sum of squares is trivially lower bounded by taking the
term corresponding to $i_{1}=\cdots=i_{q-1}=m$, which yields:
%
\begin{equation}
\label{eq:low-bnd-Am} A_{m} \ge c_{a,q,\psi}\quad\quad \mbox{with }
c_{a,q,\psi}\equiv q q! a_q^{2} \psi_0^{2q}
> 0.
\end{equation}

\textit{Step} 5: \textit{Conclusion}. In the remainder of the proof
the constants $c_{a,q,\psi,N}$ and so can change from line to
line without further mention. Plugging relation
\eqref{eq:low-bnd-Am} into \eqref{eq:low-bnd-Bn-Am} and recalling
that $N$ is a given integer whose exact value will be fixed below,
we get:
\[
\be \bigl[ B_{n}^{N-1} | \mathcal{F}_{(N-1)M} \bigr]
\ge \frac{M c_{a,q,\psi}}{n} \ge c_{a,q,\psi,N} > 0,
\]
as long as $N$ stays bounded. We then insert back this inequality
into \eqref{ineq2} and \eqref{a} in order to get:
\[
\bp \bigl( B_{n}^{N-1} \leq x | \mathcal{F}_{(N-1)M}
\bigr) \le 1 +\frac{p c_{a,q,\psi,N}}{2N} \int_0^1
x^{\afrac{1}{2(q-1)}-\afrac
{p}{2N} -1} \,\mathrm{d}x = c_{a,q,\psi,N,p} < \infty,
\]
where we have chosen $N$ such that $\frac{p}{2N} < \frac{1}{2(q-1)}$.
Iterating this bound into \eqref{ineq1}, we have thus obtained:
\[
\be \bigl[(B_{n})^{-\nfrac{p}{2}} \bigr] \le c_{a,q,\psi,N,p}^{N},
\]
which is a finite quantity. Finally, recall from Step 1 that $\be
[(B_{n})^{-\nfrac{p}{2}}]\ge\be[\|DV_{n}^{d,q}\|_{\fh}^{-p}]$, which
finishes the proof.


\section*{Acknowledgements}
We would like to thank an anonymous referee for carefully reading this
manuscript and making helpful remarks.
The research of Y. Hu is partially supported by a grant from the Simons
Foundation \#209206 and by a General Research Fund of University of
Kansas. The research of D. Nualart is supported by the NSF grant
DMS1208625. The research of F. Xu (corresponding author) is supported
by the 111 Project (B14019) and Shanghai Pujiang Program (14PJ1403300). This project
has been carried out while S. Tindel was on sabbatical at the
University of Kansas. He wishes to express his gratitude to this
institution for its warm hospitality.


%

\printhistory
\end{document}